\documentclass[12pt,a4paper,oneside]{article}
\usepackage{amsmath,amsfonts,amssymb,amsthm}

\newcommand{\Rset}{\mathbb{R}}
  \newcommand{\PP}{\mathbb{P} }
\newcommand{\Nset}{\mathbb{N} }
  \newcommand{\EE}{\mathbb{E} }
  \newcommand{\Zset}{\mathbb{Z} }
  \newcommand{\Cset}{\mathbb{C} }
  \newcommand{\e}{\textrm{e}}
        \newtheorem{thm}{Theorem}
        \newtheorem{lem}{Lemma}
        
        \newtheorem{prop}[thm]{Proposition}
        \newtheorem{hyp}{Hypotesis}[thm]
   \newcommand{\DS}{\displaystyle}

  \title{On a class of stochastic semilinear PDE's}
\author{Luigi Manca\\
        Scuola Normale Superiore di Pisa}
\begin{document}
\maketitle

\begin{center}
\textbf{ABSTRACT}
\end{center}

We consider stochastic semilinear partial differential equations with Lipschitz nonlinear terms.
 We prove existence and uniqueness of an invariant measure and the existence of a solution for the corresponding Kolmogorov equation in the space $L^2(H;\nu)$, where $\nu$ is the invariant measure. 
 We also prove the closability of the derivative operator and an integration by parts formula.
 Finally, under boundness conditions on the nonlinear term, we prove a Poincar\'e inequality, a logarithmic Sobolev inequality and the ipercontractivity of the transition semigroup. 
  



\section{Introduction and setting of the problem}        
We are concerned with the following semilinear equation perturbed by noise in the Hilbert space 
$H$ of all $2\pi$-periodic real functions 
\begin{equation} \label{1.0}
 \begin{cases}
    dX = (D_\xi^2X - X + D_\xi F(X))dt + dW, \\
    X(0)(\xi) = x(\xi), \quad \xi\in [0,2\pi],
 \end{cases}
\end{equation}
where $x\in H$, $F\in C^1(H;H)$ with $DF \in C_b(H;\mathcal{L}(H))$ and $W$ is a cylindrical Wiener process defined on a probability space 
$(\Omega, \mathcal{F}, \PP)$ with values in $H$. 
We shall denote by $\langle \cdot,\cdot \rangle$ the inner product in $H$, defined by
\begin{displaymath}
 \langle x, y \rangle = \int_0^{2\pi} x(\xi)y(\xi) d\xi, \quad x,y\in H
\end{displaymath}
and by $|\cdot|_2$ the corresponding norm.
We shall prove that \eqref{1.0} admits a unique mild solution in the space $C_W([0,T];H)$,
 consisting of all stochastic processes $X(\cdot,x) \in C([0,T];L^2(\Omega;H))$ which are adapted to $W(t)$.
We recall that a treatment of the Cauchy problem for an extensive class of Burgers-type equations can be found in \cite{GY98}.
We shall also prove the differentiability of $X(t,x)$ with respect to $x$ 
and some approximation theorems both for $X(t,x)$ and its derivative.
Through the mild solution $X(t,x)$ of \eqref{1.0} we shall define the transition semigroup $\{P_t\}_{t\geq 0}$ as
\begin{displaymath}
  P_t \varphi (x) = \EE[\varphi(X(t,x))],
\end{displaymath} 
where $\varphi:H \in \Rset$ is Borel and bounded. 
We shall prove strong Feller and irreducibility properties of the transition semigroup $\{P_t\}_{t\geq 0}$
 in order to ensure, thanks to the Doob theorem, the uniqueness of an invariant measure for the transition semigroup.
 We recall that a Borel probability measure $\nu$ is invariant for the semigroup $P_t$ if we have
\begin{equation} \label{invariant}
 \int_H P_t\varphi d\nu = \int_H \varphi d\nu
\end{equation}
for all $\varphi:H\to \Rset$ continuous and bounded.
 Then we shall present some sufficient conditions on $F$ that imply the existence (and consequently, by the Doob theorem, the uniqueness) of an invariant measure.
The existence of an invariant measure $\nu$ allow us to extend uniquely $P_t$ to a strongly continuous semigroup (still denoted by $P_t$) in $L^2(H;\nu)$.
 We shall denote by $K_2$ its infinitesimal generator.
 Then we shall show that  $K_2$ is the closure of the following differential operator 
\begin{displaymath}
  K_0\varphi(x) = \frac{1}{2} Tr[D^2 \varphi(x)] + 
  \langle (D_\xi^2 - I)x + D_\xi F(x), D\varphi(x) \rangle 
\end{displaymath}
where $Tr$ denote the trace, $D$ denote the derivative with respect to $x$ and $\varphi$ belong to a suitable subspace of $L^2(H;\nu)$ that will be rigorously defined in the following. 
This kind of result was proved for a Burgers equation with coloured noise (see \cite{DPD}). In the present situation (Lipschitz nonlinearities and a white noise perturbation) the result seems to be new.
An extensive survey on second order partial differential operators in Hilbert spaces 
can be found in the monographs
\cite{CER}, \cite{DP03}, \cite{DPZ02}.
A second new result of this paper  is the closability of the operator $D$ in $L^2(H;\nu)$ and 
 that $D(K_2)$ is included in the Sobolev space $W^{1,2}(H;\nu)$. 
This implies the integration by parts formula
\begin{equation} \label{IBPF}
 \int_H \varphi K_2 \varphi d\nu =  - \frac{1}{2} \int_H |D\varphi|^2 d\nu , \quad \varphi \in D(K_2).
\end{equation}
Moreover (but only in the case $\|DF\|_0 < 2 $) we shall show that by \eqref{IBPF} it follows a Poincar\'e-type inequality, i.e.
\begin{equation} \label{PI}
   \int_H |\varphi - \overline{\varphi}|^2 d\nu \leq
     \frac{1}{2\big(1-\frac{\|DF\|_0^2}{4}\big)}  \int_H |D\varphi|_2^2  d\nu, \quad \varphi \in W^{1,2}(H,\nu).
\end{equation}
As consequence of \eqref{PI} we shall derive that the spectrum of $K_2$ in the space 
 $L_0^2(H;\nu) =\{ \varphi \in L^2(H;\nu): \int_H \varphi d\nu = 0 \} $
 is contained in the half space $\{ \lambda \in \Cset : \Re \lambda < -(1 - \| DF\|_0^2/4) \} $.
 Moreover, we shall prove a logaritmic Sobolev inequality and consequently the hypercontractivity of $P_t$.

This paper is organized as follows: 
in the next section we introduce some notations and some functional spaces that will be used in what follows.
 Section 3 is devoted in proving existence and uniqueness of a mild solution $X(t,x)$ of problem \eqref{1.0} and  to its differentiability  with respect to $x$, and in section 4  we prove some approximation theorems.
In section 5 we introduce the transition semigroup $P_t$,
 and in sections 6, 7 we discuss the strong Feller and irreducibility properties respectively.
 In section 8 we prove the existence of an invariant measure.
In section 9 we study the infinitesimal generator $K_2$ of the semigroup $P_t$ in $L^2(H;\nu)$, 
where $\nu$ is an invariant measure for $P_t$. 
Section 10  is devoted to the integration by parts formula, and section 11  to the Sobolev space $W^{1,2}(H;\nu)$, 
i.e. the domain of the clousure of $D$ in $L^2(H;\nu)$.
 Finally, the Poincar\'e inequality, the spectral gap and the logarithmic Sobolev inequality are discussed in section 12. 


 \section{Preliminaries}              
Let us write problem \eqref{1.0} in an abstract form. For this it is convenient to consider 

the complete orthogonal system $\{e_k\}_{k\in \Zset}$  in $H$ given by
\begin{displaymath}
 e_k(\xi) = \begin{cases} \frac{1}{\sqrt{2 \pi}}\cos(k\xi), \, k\geq 0,\, \xi \in [0,2\pi],\\
                          \frac{1}{\sqrt{2 \pi}}\sin(k\xi), \, k < 0,\, \xi \in [0,2\pi].
            \end{cases}
\end{displaymath}
We represent any element $x\in H$ by its Fourier series
\begin{displaymath}
  x = \sum_{k\in \Zset} x_k e_k,\quad x_k = \langle x,e_k \rangle,
\end{displaymath}
and for any $\sigma \geq 0$ we define the set
\begin{displaymath}
   H_\#^{\sigma} = \{ x\in H : |x|_{2,\sigma} < \infty \},
\end{displaymath}
where
\begin{displaymath}
  |x|_{2,\sigma} = \big(\sum_{z\in\Zset} (1+k^2)^{\sigma/2}|x_k|^2 \big)^{1/2}.
\end{displaymath}
Now, we define a linear operator $A:D(A) \to H$ by
\begin{displaymath}
  Ax(\xi) = D^2_\xi x(\xi) - x(\xi), \, \xi \in [0,2\pi ], \, D(A) = H_\#^2(0,2\pi).
\end{displaymath}
The linear operator $A$ is selfadjoint and $ Ae_k = -(1+k^2)e_k,\quad k\in \Zset$.
Clearly we have that $|(-A)^{\sigma/2}x|_2 = |x|_{2,\sigma}$ and $|(-A + I)^{1/2}x|_2 =|D_\xi x|_2 $.
The cylindrical Wiener process $W(t)$ is formally defined  by
\begin{displaymath}
  W(t) = \sum_{z\in \Zset} \beta_k(t) e_k,\quad t\geq 0,
\end{displaymath}
where $\{\beta_k \} $ is a sequence of mutually independent real Brownian process
 in a given probability space $(\Omega, \mathcal{F},\PP)$.
Finally \eqref{1.0} can be written as              
\begin{equation} \label{1.1}
   \begin{cases}
    \DS  dX(t) = (AX(t) + D_\xi F(X(t)))dt + dW(t), \cr               
    \DS  X(0) = x \in H
  \end{cases}
\end{equation}       
In the following we will denote by $\|\cdot \|_0$ the supremum norm in the space $C(H;\mathcal{L}(H))$. Clearly the conditions on $F$ implies $\|DF\|_0 < \infty$.
We write \eqref{1.1} in the following mild form  
  \begin{equation} \label{2.1}
  X(t) = \e^{tA}x + \int_0^t D_\xi \e^{(t-s)A} F(X(t))ds + W_A(t)
  \end{equation}
 where $W_A(t)$ is the stochastic convolution
   \begin{displaymath}
 \DS    W_A(t) = \int_0^t \e^{(t-s)A} dW(s) = \sum_{k\in\Zset}\int_0^t \e^{-(t-s)k^2}e_k d\beta_k(s).
   \end{displaymath}
Notice that for any $\sigma \in [0,1/2)$ we have that $W_A(t) \in L^2(\Omega;H^\sigma_\#)$, since
 \begin{equation*}
 \DS  \| W_A(t) \|_{L^2(\Omega;H^\sigma_\#)}^2 =
            \EE |W_A(t)|_{2,\sigma} ^2 \leq \sum_{k\in\Zset} \frac{(1+k^2)^\sigma}{2(1+k^2)} < \infty.
 \end{equation*}
In order to give a precise meaning to equation \eqref{2.1}, we introduce, for any $t>0$,
the linear mapping
\begin{displaymath}
 K(t):H \to H, \, x\mapsto K(t)x,\quad K(t)x = D_\xi \e^{tA} x.
\end{displaymath}
We have
\begin{lem}
  $K(t)$ is a linear bounded mapping from $H$ into itself.
 Moreover there exists $\kappa > 0$ such that 
 \begin{equation} \label{2.3}
   | K(t)x |_2 \leq \kappa \e^{-t} t^{-1/2} | x |_2,\quad x\in H
  \end{equation}
\end{lem}
\begin{proof}
For any $t>0$ we have
  \begin{equation*}
    \DS  D_\xi \e^{tA}x = \sum_{k\in\Zset} k \e^{-(1+k^2)t} x_k e_{-k}.
  \end{equation*}
Then 
  \begin{equation*}
   \DS   |D_\xi \e^{tA}x|_2^2 =  \sum_{k\in\Zset}  k^2 \e^{-2(1+k^2)t} | x_k|^2
  \leq  \sup_{k\in \Zset} k^2 \e^{-(1+k^2)t} | x|_2^2.
 \end{equation*}
Since, as it can be easily seen,
 \begin{equation*}
   \DS   \sup_{k\in \Zset} k^2 \e^{-2(1+k^2)t} \leq \frac{1}{4\sqrt{\e}}t^{-1}\e^{-2t}
  \end{equation*}
the conclusion follows.
\qed \end{proof}
In the following will be useful the next
%
\begin{lem} \label{Henry}
 Suppose $b\geq 0$, $ \beta >0$ and that $a(t)$ is a nonnegative function
 locally integrable on $0\leq t <T$ fulfilling
\begin{displaymath}
  u(t) \leq a(t) + b \int_0^t (t-s)^{\beta - 1} u(s) ds,\quad t\in [0,T].
\end{displaymath}
  Then we have
\begin{displaymath}
 u(t) \leq a(t) + \theta \int_0^t E'_\beta (\theta(t-s))a(s)ds, \quad 0
 \leq t < T
\end{displaymath}
where
\begin{displaymath}
   \theta = (b\Gamma(\beta))^{1/\beta}, \quad E_\beta(z) =
   \sum_{n=0}^\infty \frac{z^{n\beta}}{\Gamma(n\beta + 1)}, \quad
   E'_\beta(z)=\frac{d}{dz} E_\beta(z).
\end{displaymath}
Moreover
\begin{displaymath}
  E'_\beta(z) \thicksim \frac{z^{\beta - 1}}{\Gamma(\beta)} \;\text{as $z\to
  0^+$}, \quad E'_\beta(z)\thicksim E_\beta(z) \thicksim \frac{\e^z}{\beta} \;
  \text{as $z\to + \infty$},
\end{displaymath}
and if $a(t)=a$, constant, then $u(t)\leq a E_\beta(\theta t)$.
\end{lem}
\begin{proof}
See e.g. Lemma 7.1.1 on \cite{Hen}.
\qed \end{proof}
%
%
\section{The mild solution $X(t,x)$ and its differentiability}
\begin{thm}
For any $x \in H$ and $T>0$ there exists a unique mild solution
$X\in C_W ([0,T];H)$ of equation \eqref{1.1}.
\end{thm}
\begin{proof}
Existence and uniqueness of a solution of equation \eqref{1.1} follows easily by the fixed point method in the space $C_W ([0,T];H) $.
\qed \end{proof}
%
%
%
%
We prove here that the mild solution $X(t,x)$ of \eqref{2.1}
is differentiable with respect to $x$ and that for any $h\in H$ it
holds
\begin{displaymath}
  DX(t,x)\cdot h=\eta^h(t,x),
\end{displaymath}
where $\eta^h(t,x)$ is the mild solution of the equation
\begin{equation}\label{ED.1}
 \begin{cases}
  \DS \frac{d}{dt}\eta^h(t,x) = A\eta^h(t,x) + D_\xi( DF(X(t,x)\cdot\eta^h(t,x)) \cr
     \eta^h(0,x)=h
 \end{cases}
\end{equation}
This means that $\eta^h(t,x)$ is the solution of the integral
equation 
\begin{equation} \label{ED.2}
  \eta^h(t,x) =
    \e^{tA}h + \int_0^t K(t-s) DF(X(s,x))\cdot \eta^h(s,x) ds,
    \quad t\geq 0.
\end{equation} 
\begin{thm} \label{ED} 
 Assume that $X(t,x)$ is the solution of equation \eqref{2.1}. Then it is differentiable with
 respect to $x$ $ \PP$-a.s., and for any $h\in H$ we have 
 \begin{equation} \label{ED.3}
  DX(t,x)\cdot h = \eta^h(t,x),\quad \PP-\text{a.s.}
 \end{equation}
 and 
 \begin{equation} \label{ED.4}
  |\eta^h(t,x)|_2 \leq \e^{\big(\frac{\|DF\|_0^2}{4} -1\big)t} |h|_2, \quad  t\geq 0
 \end{equation}
\end{thm}
\begin{proof}
 Arguing as in the proof of Theorem 1, we notice that \eqref{ED.1} has a unique mild solution
 $\eta^h(t,x)$ in $C_W([0,T];H)$.
  Let us prove \eqref{ED.4}.
  By multiplying both sides of \eqref{ED.1} by $\eta^h(t,x)$ and integrating on $[0,2\pi]$ we have
\begin{displaymath}
 \frac{1}{2}\frac{d}{dt}|\eta^h(t,x)|_2^2 =
  \langle A\eta^h(t,x),\eta^h(t,x) \rangle + \langle D_\xi\big(
  DF(X(t,x))\cdot \eta^h(t,x)\big),\eta^h(t,x) \rangle.
\end{displaymath}
Integrating by parts and applying the H\"older inequality we find
\begin{multline*}
 \frac{1}{2}\frac{d}{dt}|\eta^h(t,x)|_2^2 \leq \\
 \leq  \langle A\eta^h(t,x),\eta^h(t,x) \rangle
    + \frac{\| DF \|_0^2}{4}|\eta^h(t,x)|_2^2 + |D_\xi\eta^h(t,x)|_2^2 = \\
 = \big( \frac{\| DF \|_0^2}{4} -1 \big) |\eta^h(t,x)|_2^2.
 \end{multline*}
 Then \eqref{ED.4} follows by Gronwall's lemma.

 Now we prove that $\eta^h(t,x)$ fulfills \eqref{ED.3}.
 For this fix $T >0$, $x,h \in H$ such that $|h|_2 \leq 1$.
 We claim that there exist a constant $C_T > 0$ and a function
 $\sigma_T(\cdot):H\to\Rset^+$, with $\sigma_T(h)\to 0$ as $h\to 0$,
 such that
\begin{displaymath}
  |X(t,x+h) - X(t,x) - \eta^h(t,x)|_2 \leq  C_T\sigma_T(h)|h|_2, \quad \PP-\text{a.s.}.
\end{displaymath}
 Setting
\begin{displaymath}
   r_h(t,x)=X(t,x+h) - X(t,x) - \eta^h(t,x),
\end{displaymath}
    $r_h(t,x)$ satisfies the equation
\begin{multline*}
   r_h(t,x) = \int_0^t K(t-s)\big[F(X(s,x+h)) - F(X(s,x))]ds + \\
   - \int_0^t K(t-s)DF(X(s,x))\cdot \eta^h(s,x) ds.
\end{multline*}
   Consequently we have that
\begin{eqnarray*}
 r_h(t,x) &=& \int_0^t K(t-s)\int_0^1 DF(\rho (\zeta,s))d\zeta \big(X(s,x+h) + X(s,x)\big)ds + \\
 & &- \int_0^t K(t-s)DF(X(s,x))\cdot \eta^h(s,x) ds = \\
 &=& \int_0^t K(t-s)\int_0^1 DF(\rho (\zeta,s))d\zeta \cdot r_h(s,x)ds + \\
 & &+ \int_0^t K(t-s) \int_0^1 \big( DF(\rho (\zeta,s)) 
  -DF(X(s,x))\big)d\zeta \cdot \eta^h(s,x) ds.
\end{eqnarray*}
 where $\rho (\zeta,s) = \zeta X(s,x+h) + (1-\zeta)X(s,x)$.
 Notice that since $F\in C^1_b(H)$ and $X(t,x)$ is continuous with
 respect to $x$ uniformly in $[0,T]$, there exists a function
 $\sigma_T:H \to \Rset^+$ such that $\sigma_T \to 0$ as $h\to 0$ and
\begin{equation}\label{ED.5}
 |DF(\rho(\zeta,s)) - DF(X(s,x))|_2 \leq \sigma_T(h).
\end{equation}
Setting
\begin{displaymath}
   \gamma_T=\sup_{t\in[0,T]}\e^{(\|DF\|_0^2/4-1)t},
\end{displaymath}
  and taking into account \eqref{ED.4},\eqref{ED.5}, we find
\begin{multline*}
 |\int_0^t K(t-s)\int_0^1 \big(DF(\rho(\zeta,s)) - DF(X(s,x))\big)d\xi
 \cdot \eta^h(s,x)ds|_2 \leq \\
 \leq \kappa \int_0^t \e^{-(t-s)}(t-s)^{-1/2}ds\sigma_T(h)|h|_2  \leq
 \kappa\Gamma(\frac{1}{4})\gamma_T \sigma_T(h) |h|_2.
\end{multline*}
It follows that
\begin{displaymath}
   |r_h(t,x)|_2 \leq  \|DF\|_0\int_0^t \e^{-(t-s)}(t-s)^{-1/2}|r_h(s,x)|_2  ds +
   \kappa\Gamma(\frac{1}{4})\gamma_T \sigma_T(h) |h|_2,
\end{displaymath}
and thus by Lemma \ref{Henry} we have $|r_h(t,x)|_2 \leq \kappa\Gamma(1/2)\gamma_T E_{1/2} (\theta T) \sigma_T(h),
 |h|_2$,
 where $\theta = \big(\|DF\|_0 \Gamma(1/2)\big)^2$. This implies \eqref{ED.3}.
 \qed \end{proof}
%
%
%
%
%
%
%
\section{Approximation of $X(t,x)$ and $\eta^h(t,x)$}
In this section we consider the approximated problem 
\begin{equation} \label{approx.1}
   \begin{cases}
    \DS  dX_n(t) = (AX_n(t) + D_{\xi,n} (F(X_n(t)))dt + dW(t), \cr
    \DS  X_n(0) = x \in H,
  \end{cases}
\end{equation}
where $D_{\xi,n} \in \mathcal{L}(H)$ is defined by $D_{\xi,n} = D_\xi \circ P_n$
and $P_n$ is the projection of $H$ into the linear span of $\{e_{-n},\ldots, e_n \}$. 
We also consider problem \eqref{approx.1} in its mild form, i.e. 
\begin{equation} \label{approx.2}
  X_n(t) = \e^{tA}x + \int_0^t K_n(t-s)F(X_n(s))ds + W_A(t),
\end{equation}
where $K_n(t) = D_{\xi,n} \e^{tA}$. 
Notice that $D_{\xi,n}\circ F:H \to H$ is a nonlinear Lipschitz continuos function, 
and so, as it is well know (see, for example, \cite{DPZ}), problem \eqref{approx.1} admits a mild solution in $C_W([0,T];H)$.
Moreover, for any $n\in \Nset$, 
 $t\geq 0$ we have that  $K_n(t) \in \mathcal{L}(H)$  and it holds 
\begin{equation} \label{approx.3}
 \begin{array}{l}
    \DS \|K_n(t)\|_{\mathcal{L}(H)} < \|K(t)\|_{\mathcal{L}(H)}, \cr
    \DS K_n(\cdot) \to K(\cdot) \text{ in } C([t_0,T];\mathcal{L}(H)),\, 0<t_0<T.
 \end{array}    
\end{equation}
%
%
We have 
\begin{thm} \label{approx}
 If $X_n(t,x)$ and $X(t,x)$ are the solutions of problem \eqref{approx.2} and \eqref{2.1} respectively, 
 then 
\begin{equation}\label{approx.30} 
   \lim_{n\to \infty} X_n(\cdot,x) = X(\cdot,x),\quad \text{in } C_W([0,T];H).
\end{equation}
\end{thm}
\begin{proof}
Let us fix $\varepsilon > 0$. Taking into account \eqref{2.1}, \eqref{approx.2} we have 
\begin{multline*}
 X(t,x) - X_n(t,x) = \int_0^t (K(t-s) - K_n(t-s)) F(X(s))ds + \\
        + \int_0^t K_n(t-s)(F(X(s)) - F(X_n(s))ds  
\end{multline*}
Moreover, taking into account \eqref{approx.3}, for all $0<t_0< t\leq T$ it holds 
\begin{multline*}
 |\int_0^t (K(t-s) - K_n(t-s)) F(X(s))ds|_2  \leq\\
   \leq  \int_0^t \|K_n(s) - K(s)\|_{\mathcal{L}(H)} ds \|DF\|_0 \sup_{0\leq t \leq T} |X(t)|_2 \leq         \\
 \leq \big(2\int_0^{t_0}\|K(s)\|_{\mathcal{L}(H)} + \int_{t_0}^t \|K(s) - K_n(s)\|_{\mathcal{L}(H)}ds\big) 
  \|DF\|_0 \sup_{0\leq t \leq T}|X(t)|_2         \leq \\
\leq \big( 4 \kappa \sqrt{t_0} + T \sup_{t_0 \leq t \leq T}\|K(t) - K_n(t)\|_{\mathcal{L}(H)} \big) \|DF\|_0 \sup_{0\leq t \leq T}|X(t)|_2 
\end{multline*}
and
\begin{multline*}
  |\int_0^t K_n(t-s)(F(X(s)) - F(X_n(s))ds|_2 \leq \\
 \leq   \kappa \|DF\|_0 \int_0^t (t-s)^{-1/2}\e^{-(t-s)}|X(s) - X_n(s)|_2ds.
\end{multline*} 
Then, by Lemma \ref{Henry}, it follows 
\begin{multline*}
 |X(t) - X_n(t)|_2 \leq 
  \big( 4 \kappa \|DF\|_0 \sup_{0\leq t \leq T} |X(t,x)|_2 \sqrt{t_{0}} + \\
  + T\|DF\|_0 \sup_{0\leq t \leq T}|X(t,x)|_2 
    \|K(\cdot) -  K_n(\cdot)\|_{C([t_0,T],\mathcal{L}(H))}\big)E_{1/2}(\theta T),
 \end{multline*}
and consequently
\begin{multline*}
\sup_{t\in [0,T]}\EE|X(t) - X_n(t)|_2^2 
       \leq 2 \big( 16 \kappa^2 \|DF\|_0^2  t_{0} + \\
      + T^2\|DF\|_0^2 
      \|K(\cdot) -  K_n(\cdot)\|_{C([t_0,T],\mathcal{L}(H))}^2\big)E_{1/2}(\theta T)^2\|X(\cdot,x)\|_{C_W([0,T];H)}
\end{multline*}                          
Now taking 
\begin{displaymath}
  t_0 <  \big( 16\kappa^2\|DF\|_0^2  E_{1/2}(\theta T)^2\|X(\cdot,x)\|_{C_W([0,T];H)}\big)^{-1} \frac{\varepsilon}{4}
\end{displaymath}
and $n$ such that 
\begin{displaymath}
 \sup_{t\in [t_0,T]}  \|K(\cdot) -  K_n(\cdot)\| <  
  \big( T^2\|DF\|_0^2  E_{1/2}(\theta T)^2\|X(\cdot,x)\|_{C_W([0,T];H)}\big)^{-1/2} (\frac{\varepsilon}{4})^{1/2}
\end{displaymath}   
we find 
\begin{displaymath}
  \sup_{t\in [0,T]}\EE|X(t) - X_n(t)|_2^2 < \varepsilon.
\end{displaymath}  
Theorem \ref{approx} is proved.
\qed \end{proof}
%
%
Denote with $\eta_n^h(t,x)$ the mild solution of problem 
\begin{equation}\label{approx.4}
 \begin{cases}
  \DS \frac{d}{dt}\eta_n^h(t,x) = A\eta_n^h(t,x) + D_{\xi,n}( DF(X_n(t,x)\cdot\eta_n^h(t,x)), \cr
     \eta_n^h(0,x)=h.
 \end{cases}
\end{equation}
It is well know that the solution $X_n(t,x)$ of problem \eqref{approx.1} it is differentiable with respect to $x$ $\PP-$a.s. (see, for example, \cite{DPZ}), and that 
\begin{displaymath}
   \langle DX_n(t,x),h\rangle = \eta_n^h(t,x), \quad h\in H, t\geq 0.
\end{displaymath}
Moreover it is easy to see that \eqref{ED.4} still holds for $\eta_n^h(t,x)$. We have also the next
\begin{thm}
 If $\eta^h(t,x)$ and $\eta_n^h(t,x)$ are the solutions of problems \eqref{ED.1}, \eqref{approx.4} respectively, then for all $h\in H$ 
\begin{equation} \label{approx.5}
  \lim_{n\to\infty} \eta_n^h(t,x) = \eta^h(t,x) 
\end{equation} 
in  $ C_W([0,T];H) $
\end{thm} 
\begin{proof}
 The proof is similar to that of Theorem \ref{approx}
 \qed \end{proof}
%
%
%
\section{The transition semigroup}
The transition semigroup corresponding to the mild solution $X(t,x)$ of \eqref{2.1} is defined by 
\begin{equation} \label{TS.1}
    P_t\varphi(x)= \EE[\varphi(X(t,x))],\quad \varphi \in B_b(H),\, t\geq 0,\, x\in H.
\end{equation}
Let us also consider the approximating semigroup 
\begin{equation} \label{TS.2}
      P_t^n\varphi(x)= \EE[\varphi(X_n(t,x))],\quad \varphi \in B_b(H),\, t\geq 0,\, x\in H,
\end{equation}
for all $n\in \Nset$, where $X_n(t,x)$ is the solution of \eqref{approx.2}.
We have obviously
\begin{displaymath}
   \|P_t\varphi\|_0 \leq \|\varphi\|_0, \quad \varphi \in B_b(H),
\end{displaymath}
and by the dominated convergence theorem it follows that
\begin{displaymath}
  \lim_{n \to \infty} P_t^n\varphi(x) = P_t\varphi(x),\quad \varphi \in C_b(H),\, x\in H.
\end{displaymath}
If $F\in C_b^1(H;H)$, by Theorem \ref{ED} we have that, for all $\varphi \in C_b^1(H)$, $P_t\varphi(x)$ and $P_t^n\varphi(x)$ are differentiable with respect to $x$ and it holds
\begin{displaymath}
   \langle DP_t\varphi(x), h\rangle = \EE \langle D\varphi(X(t,x)), \eta^h(t,x)\rangle ,\quad h\in H, 
\end{displaymath}
\begin{displaymath}
   \langle DP_t^n\varphi(x), h\rangle = \EE \langle D\varphi(X_n(t,x)), \eta^h_n(t,x)\rangle, \quad h\in H.
\end{displaymath}
Moreover, by Theorem \ref{approx} and \eqref{approx.5} it follows that for all $\varphi \in C_b^1(H)$, $h\in H$, 
\begin{displaymath}
     \lim_{n \to \infty}  \langle DP_t^n\varphi(x), h\rangle =  \langle DP_t\varphi(x), h\rangle
\end{displaymath} 
in $C([0,T]; \Rset)$.
%
%
%
\section{Strong Feller property}
In order to prove the strong Feller property of the transition semigroup $P_t$,
i.e for all $\varphi \in B_b(H)$, $t>0$,  it follows that $P_t\varphi \in C_b(H)$,
 we shall use the Bismut-Elworthy formula (see \cite{DPEZ}). 
 Since $D_\xi F$ is not Lipschitz continuous, we will apply the Bismut-Elworthy formula to the approximated transition semigroup $P_t^n$, defined in \eqref{TS.2}, and then te shall let $n \to \infty$.
\begin{lem}
 If $\varphi \in C_b^2(H)$ and $t>0$  we have, for all $n\in \Nset$, $P_t^n\varphi \in
 C_b^1(H)$ and, for any $h \in H$,  
\begin{equation}\label{3.3}
    \langle DP_t^n\varphi(x), h \rangle
  = \frac{1}{t}\EE \Big[ \varphi(X_n(t,x)) \int_0^t \langle \eta^h_n(s,x), dW(s)\rangle\Big].
\end{equation}
\end{lem}
\begin{proof} See \cite{DPEZ}.
\qed \end{proof}
Formula \eqref{3.3} remains true also for $\varphi \in C_b(H)$, since
  we can pointwise approximate a $C_b(H)$-function by  a sequence of $C_b^2(H)$-functions. 
\begin{thm} \label{SF}
 The transition semigroup $P_t$ defined in \eqref{TS.1} is strong Feller.
\end{thm}
\begin{proof}
Step 1. If $\varphi \in C_b^2(H)$, for all $t>0$ we have
\begin{displaymath}
|DP_t\varphi(x)|_2 \leq  t^{-1}\frac{\sqrt{2}}{\|DF\|_0}(\e^{\frac{\|DF\|_0^2}{2}t} -1)^{1/2} \| \varphi  \|_0.
\end{displaymath}
In fact by \eqref{3.3}, using the H\"older inequality and recalling
\eqref{ED.4}, for all $n\in \Nset$ we have
\begin{displaymath}
 |\langle DP_t^n\varphi(x), h \rangle |^2
      \leq t^{-2} \| \varphi \|_0^2 \int_0^t |\eta^h_n(s,x)|_2^2 ds \leq
\end{displaymath}
\begin{displaymath}
   \leq t^{-2} \| \varphi \|_0^2 \int_0^t \e^{\frac{\|DF\|_0^2}{2}s}|h|_2^2ds =
  t^{-2} \| \varphi \|_0^2 \frac{2}{\|DF\|_0^2}( \e^{\frac{\|DF\|_0^2}{2}t} - 1)|h|_2^2.
\end{displaymath}
Now, letting $n\to \infty$, the conclusion holds for the arbitrariness of $h$.\\
Step 2. For any $\varphi \in B_b(H)$, $t>0$ and $x,y \in H$  it holds  
\begin{equation} \label{3.4}
|P_t \varphi(x) - P_t \varphi(y)| \leq t^{-1} \frac{\sqrt{2}}{\|DF\|_0}( \e^{\frac{\|DF\|_0}{2} t} - 1)^{1/2} \| \varphi \|_0|x-y|_2
\end{equation}
In order to prove the step we need to approximate $\varphi$ by a
sequence of $C_b^2(H)$-functions. Since  $C_b^2(H)$ is not dense in
$B_b(H)$, we will use a suitable pointwise approximation. Fix $t>0$
and $x,y \in H$. Let us define a signed measure $\zeta$ setting
$\zeta = \lambda_{t,x} - \lambda_{t,y}$, where $\lambda_{t,x},
\lambda_{t,y}$ are the law of $X(t,x)$ and $X(t,y)$ respectively, and
consider a sequence $\{\varphi_n \} $ of $C_b^2(H)$-functions such
that
\begin{displaymath}
 \lim_{n\to \infty} \varphi_n(x) = \varphi(x) \quad\text{$\zeta$-a.s.},\quad
\| \varphi_n \|_0 \leq \| \varphi \|_0 \quad \forall n\in \Nset .
\end{displaymath}
By step 1 we have
\begin{displaymath}
\DS |P_t \varphi_n(x) - P_t \varphi_n(y)|
    \leq \sup_{0\leq \theta \leq 1}
     \| DP_t\varphi_n (\theta x + (1-\theta)y))\|_{\mathcal{L}(H)}|x-y|_2 \leq
\end{displaymath}
\begin{displaymath}
   \leq  t^{-1} \frac{\sqrt{2}}{\|DF\|_0}( \e^{\frac{\|DF\|_0}{2} t} - 1)^{1/2} \|\varphi_n \|_0 |x-y|_2.
\end{displaymath}
By the dominate convergence theorem, it follows that \eqref{3.4} holds
and so $P_t\varphi \in C_b(H)$ as claimed. Theorem \ref{SF} is proved.
\qed \end{proof}
%
\section{Irreducibility}
A basic tool for proving irreducibility of $P_t$ is the approximate
controllability of the following controlled system  
\begin{equation} \label{4.1}
  \begin{cases}
    y'(t)=Ay(t) + D_\xi F(y(t)) + u(t) \cr
    y(0) = x
  \end{cases}
\end{equation}
where $u \in L^2([0,T];H)$.
Let us denote by $y(\cdot,x;u)$ the mild solution of \eqref{4.1}, that is the solution of the
integral equation  
\begin{equation} \label{4.2}
    y(t)=\e^{tA}x + \int_0^t K(t-s)F(y(s))ds + \sigma_u(t),
\end{equation}
where
\begin{displaymath}
\sigma_u(t)=\int_0^t \e^{(t-s)A}u(s)ds.
\end{displaymath}
We say that the sistem \eqref{4.1} is approximatively controllable
if for any $\varepsilon >0,\,T>0$, $x,\, z \in H$, there exists $u
\in L^2([0,T];H)$ such that 
\begin{equation} \label{4.3}
    |y(T,x;u)-z|\leq \varepsilon.
\end{equation}
We have
\begin{lem}
The system \eqref{4.1} is approximatively controllable.
\end{lem}
\begin{proof}
Let be $\varepsilon >0,\, T>0, \, x,\,z\in H$. we have to show that there exists $u \in L^2([0,T];H)$
such that \eqref{4.3} holds.\\
Step 1. The mapping
\begin{displaymath}
\sigma: L^2([0,T];H) \to C_0([0,T];H) \quad  u\mapsto \sigma_u,
\end{displaymath}
where
\begin{displaymath}
 C_0([0,T];H) = \{x\in C([0,T];H):x(0)=0 \}
\end{displaymath}
has dense range. In fact is easy to check that the set
\begin{displaymath}
D_0 = \{\varphi \in  C^1([0,T];D(A)): x(0) = 0 \}
\end{displaymath}
is dense in $ C_0([0,T];H)$. Now let $\varphi \in D_0$ and set
\begin{displaymath}
u(t) = \varphi(t) - A\varphi(t) - D_\xi\varphi(t).
\end{displaymath}
It is clear that $\sigma_u = \varphi$, so the range of $\sigma$ is dense as claimed.\\
Step 2. Conclusion.\\
Choose $\psi \in C([0,T];H)$ such that $\psi(0)=0,\, \psi(T)=z$, for istance
\begin{displaymath}
\psi(t)= \frac{T-t}{T} x + \frac{t}{T}z,\, t\in [0,T],
\end{displaymath}
and set
\begin{displaymath}
   g(t)= \psi(t)-\e^{tA}x-\int_0^t K(t-s)F(\psi(s))ds,\, t\in [0,T].
\end{displaymath}
Now, given $\varepsilon > 0$, by Step 1 there exists $u\in  L^2([0,T];H)$ such that
\begin{displaymath}
  |\sigma_u(t) - g(t)|\leq C,\, t \in [0,T],
\end{displaymath}
where the constant $C$ will be choosen later. Let us show that
\eqref{4.3} holds.
 In fact, let $y(\cdot,x;u)$ be the solution of \eqref{4.1}.
 By \eqref{2.3} we have
\begin{displaymath}
  |y(t)-\psi(t)|_2 \leq \int_0^t |K(t-s)\big(F(y(s))- F(\psi(s))\big)|_2 ds + |\sigma_u(t) - g(t)|_2 \leq
\end{displaymath}
\begin{displaymath}
 \leq \kappa \|DF\|_0 \int_0^t \e^{-(t-s)} (t-s)^{-1/2}|y(s)-\psi(s)|_2 ds + |\sigma_u(t) - g(t)|_2.
\end{displaymath}
 Then by Lemma \ref{Henry} it follows that
\begin{displaymath}
  |y(t)-\psi(t)|_2 \leq C E_{1/2}(\theta t)
\end{displaymath}
 and consequently
\begin{displaymath}
  |y(T)- z|_2 \leq C E_{1/2}(\theta T).
\end{displaymath}
Now it is enough to choose $C <  E_{1/4}(\theta T)^{-1}
\varepsilon$.
\qed \end{proof}
%
%
%
%
%
\begin{thm} \label{IRR}
   The transition semigroup $P_t$ defined in \eqref{TS.1}  is irriducible.
\end{thm}
\begin{proof}
Let be $\varepsilon,\,T>0, \, x,\,z \in H$.
 We have to show that  
\begin{equation} \label{4.4}
   P_t \chi_{B^c(z,\varepsilon)}(x) = \PP(|X(t,x)-z|_2>\varepsilon)<1,
\end{equation}
where $X(t,x)$ is the solution of \eqref{2.1}.
 For this purpose we choose a control $u\in L^2([0,T];H)$ such that $|y(T,x;u)-z|_2 \leq \varepsilon/2$,
where $y$ is the solution of \eqref{4.2}.
Since
\begin{displaymath}
  |X(T,x)-z|_2 \leq |X(T,x)-y(T,x)|_2 + \frac{\varepsilon}{2},
\end{displaymath}
we have  
\begin{equation} \label{4.5}
   \PP(|X(T,x)-z|_2>1) \leq \PP(|X(T,x)-y(T,x)|_2> \frac{\varepsilon}{2}).
\end{equation}
But by \eqref{2.3} it holds
 \begin{multline*}
  |X(t,x)-y(t)|_2 \leq \\
  \leq \int_0^t |K(t-s)\big(F(X(s,x))-F(y(s,x))\big)|_2ds +
                                  |W_A(t) - \sigma_u(t)|_2 \leq \\
  \leq \kappa \|DF\|_0 \int_0^t \e^{-(t-s)}(t-s)^{-1/2} |X(s,x)-y(s,x)|_2ds +
                                  |W_A(t) -  \sigma_u(t)|_2.
 \end{multline*}
consequently, by Lemma \ref{Henry}, it follows that
\begin{displaymath}
  |X(t,x)-y(t)|_2 \leq
        |W_A(t) - \sigma_u(t)|_2 +
        \theta \int_0^t E'_{1/2}(\theta (t-s))|W_A(s) - \sigma_u(s)|_2.
\end{displaymath}
Moreover, since $W_A(\cdot)$ is a nondegenerate continuous Gaussian random variable, we have that $\PP(\sup_{t\in[0,T]}|W_A(t) - \sigma(t)|_2 > \varepsilon)$ $< 1$. This implies that
\begin{multline*}
\PP(|X(T,x)-y(T)|_2 > \varepsilon) \leq \\
 \PP( |W_A(t) - \sigma_u(t)|_2 +
        \theta \int_0^t E'_{1/2}(\theta (t-s))|W_A(s) -
    \sigma_u(s)|_2>\varepsilon) < 1,
 \end{multline*}
 and therefore \eqref{4.4} is proved.
\qed \end{proof}
%
%
\section{Existence and uniqueness of an invariant measure}
In this section we shall assume that
\begin{hyp}
\begin{align}
 \|F\|_0 &< \infty \label{Hyp1} \\
   \text{or} \notag \\
 \|DF\|_0 &< 2 \label{Hyp2} \\
   \text{or} \notag \\
   F\in C^1(\Rset;\Rset) &\text{ and } \|F'\|_0 < \infty \label{Hyp3}
\end{align}
\end{hyp}
In order to prove the existence of an invariant measure we set
$Y(t)=X(t,x)-W_A(t)$, where $X(t,x)$ is the solution of problem
\eqref{1.1}.  
Since $Y(t)$ is the solution of the integral equation
\begin{displaymath}
  Y(t) = \e^{tA}x + \int_0^t K(t-s)F(Y(s)+W_A(s))ds,
\end{displaymath}
it follows easily that $Y(t)$ is the strong solution of  
%
%
\begin{equation} \label{E.00}                                  
 \begin{cases}
      \DS \frac{d}{dt}Y(t) = AY(t) + D_\xi F(Y(t)+W_A(t)), \\
      Y(0) = x.
 \end{cases}
\end{equation}
Multiplying both sides of \eqref{E.00} by $Y(t)$ and integrating over
$[0,2\pi]$ we find     
\begin{equation}\label{E.0}
  \frac{1}{2}\frac{d}{dt}|Y(t)|_2^2 = \langle AY(t),Y(t) \rangle +
  \langle D_\xi F(Y(t)+ W_A(t)),Y(t) \rangle.
\end{equation}
We have the next 
\begin{lem} \label{lemmaE1.1}
Assume that \eqref{Hyp1} holds. Then for all $ 0\leq \varepsilon \leq 1$ it holds
\begin{equation}\label{E1.1}
 |Y(t)|_2^2 + \int_0^t \e^{-2(1-\varepsilon)(t-s)}|D_\xi Y(s)|_2^2 ds \leq
  |x|_2\e^{-2(1-\varepsilon)t} + \|F\|_0^2\int_0^t \e^{-2(1-\varepsilon)(t-s)}ds
\end{equation}
\end{lem}
\begin{proof}
Fix $0\leq \varepsilon\leq 1$. By \eqref{Hyp1} and \eqref{E.0} it holds
\begin{multline*}
   \frac{1}{2}\frac{d}{dt}|Y(t)|_2^2 \leq \langle AY(t),Y(t) \rangle
   + \|F\|_0^2 +\frac{|D_\xi Y(t)|_2^2}{2} = \\
   = -| Y(t)|_2^2 - \frac{|D_\xi Y(t)|_2^2}{2} +  \|F\|_0^2 \leq
    -(1-\varepsilon) | Y(t)|_2^2 - \frac{|D_\xi Y(t)|_2^2}{2} +  \|F\|_0^2
\end{multline*}
 Now \eqref{E1.1} follows by the Gronwall lemma.
\qed \end{proof}
\begin{lem}
 Assume that \eqref{Hyp2} holds. Then for all $ \|DF\|_0^2/4 < \varepsilon \leq 1 $ it holds 
\begin{multline} \label{E1.2}
 |Y(t)|_2^2 + (1- \frac{\|DF\|_0^2}{4\varepsilon})\int_0^t \e^{-(1-\varepsilon)(t-s)} |D_\xi Y(s)|_2^2 ds \leq \\
 \leq \e^{-(1-\varepsilon)t}|x|_2^2 + \frac{\|DF\|_0^2}{1-\frac{\|DF\|_0^2}{4\varepsilon }}\big(\int_0^t \e^{-(1-\varepsilon)(t-s)} |W_A(s)|_2^2 ds \big)
 \end{multline}
 \end{lem}
\begin{proof}
Fix  $\|DF\|_0^2/4 < \varepsilon \leq 1 $. Integrating by parts and applying the Young inequality we find, for all $M>0$, 
\begin{multline*}
 |\langle D_\xi F(Y(t) + W_A(t)),Y(t) \rangle| \leq
 |F(Y(t) + W_A(t))|_2|D_\xi Y(t)|_2 \leq \\
 \leq  \|DF\|_0|Y(t)|_2|D_\xi Y(t)|_2 + \|DF\|_0|W_A(t)|_2|D_\xi Y(t)|_2 \leq \\
 \leq \varepsilon |Y(t)|_2 + \frac{\|DF\|_0^2}{4\varepsilon }|D_\xi Y(t)|_2 + \frac{\|DF\|_0^2}{2M}|W_A(t)|_2 + \frac{M}{2}|D_\xi Y(t)|_2.
\end{multline*}
Then by \eqref{E.0} we have
\begin{displaymath}
 \frac{1}{2}\frac{d}{dt}|Y(t)|_2^2 \leq -(1-\varepsilon)|Y(t)|_2^2 + (\frac{M}{2} + \frac{\|DF\|_0^2}{4} -1)|D_\xi Y(t)|_2 +  \frac{\|DF\|_0^2}{2M}|W_A(t)|_2
\end{displaymath}
Since $\|DF\|_0 < 2 $ we can set $M = 1 - \|DF\|_0^2/4\varepsilon $ and so we find
\begin{displaymath}
 \frac{1}{2}\frac{d}{dt}|Y(t)|_2^2 \leq -(1-\varepsilon)|Y(t)|_2^2  - \frac{1}{2} (1- \frac{\|DF\|_0^2}{4\varepsilon})|D_\xi Y(t)|_2 + 
  \frac{\|DF\|_0^2}{1-\frac{\|DF\|_0^2}{4\varepsilon}} |W_A(t)|_2^2. 
\end{displaymath}
Now applying the Gronwall lemma we find \eqref{E1.2}.
\qed \end{proof}
\begin{lem}          
 If \eqref{Hyp3} holds, then $DF = F'$ and  for all $0\leq \varepsilon \leq 1$ it holds 
\begin{multline} \label{E1.3}
 |Y(t)|_2^2 + \int_0^t \e^{-2(1-\varepsilon)(t-s)}|D_\xi Y(s)|_2^2 ds \leq \\
      \leq   \e^{-2(1-\varepsilon)t}|x|_2^2 + \|DF\|_0^2 \int_0^t \e^{-2(1-\varepsilon)(t-s)}|W_A(s)|_2^2ds 
\end{multline}
\end{lem}
\begin{proof}
Since $F \in C^1(\Rset;\Rset)$, it is easy to see that for all $x \in H$, $\xi \in [0,2\pi]$, we have $F(x)(\xi)=F(x(\xi))$ and therefore $DF = F'$.
We have also that for all $x\in H^1$ it holds 
\begin{equation} \label{E1.4}
  \langle D_\xi F(x),x \rangle = 0.
\end{equation}
In fact we have
\begin{multline*}
   \langle D_\xi F(x),x \rangle = - \langle  F(x),D_\xi x \rangle = \\
 =  \int_0^{2\pi} F(x(\xi))D_\xi x(\xi) d\xi = \int_0^{2\pi} D_\xi \big(\int_0^{x(\xi)} F(\xi')d\xi'\big) d\xi = 0.
\end{multline*}
This implies that for all $x\in H^1$, $y\in H$ it holds
\begin{equation} \label{E1.5}
 |\langle D_\xi F(x+y), x \rangle| \leq \|DF\|_0 |y|_2 |D_\xi x|_2
\end{equation}
In fact, taking into account \eqref{E1.4} we have
\begin{displaymath}
 | \langle D_\xi F(x+y), x \rangle  |= | \langle  F(x+y),D_\xi x \rangle | 
 = | \langle  F(x+y) - F(x) ,D_\xi x \rangle |
\end{displaymath}
that implies \eqref{E1.5}.
Now fix $\varepsilon \geq 0$. Then,  by \eqref{E.0}, it follows 
\begin{multline*}
  \frac{1}{2}\frac{d}{dt}|Y(t)|_2^2 \leq \langle AY(t),Y(t) \rangle +
  \|DF\|_0 | W_A(t)|_2|D_\xi Y(t)|_2 =    \\
    = -|Y(t)|_2^2 -\frac{|D_\xi Y(t)|_2^2}{2} + \frac{ \|DF\|_0 }{2} | W_A(t)|_2^2 \leq  \\
  \leq  -(1-\varepsilon)|Y(t)|_2^2 -\frac{|D_\xi Y(t)|_2^2}{2} + \frac{ \|DF\|_0 }{2} | W_A(t)|_2^2,
\end{multline*}
and so applying the Gronwall lemma  yields \eqref{E1.3}. 
\qed \end{proof}
%
%
%
Now we are able to prove the main result of this section.
\begin{thm}
Let  $X(t,x)$ be the mild solution of problem \eqref{1.1}. 
If hypotesis 5.1 holds then there exists a unique invariant measure  for the transition semigroup $P_t$ defined in \eqref{TS.1}.
\end{thm}
\begin{proof}
Since for Theorem \ref{SF} and Theorem \ref{IRR} the transition semigroup $P_t$ is strong Feller and irreducible,
 it is sufficient to prove the existence of an invariant measure (see \cite{DPZ96}). 
So, fix $x\in H$ and denote with $\lambda_{t,x}$ the law of $X(t,x)$ and by $\mu_{T}$ the measure
\begin{displaymath}
  \mu_{T} = \frac{1}{T} \int_0^T \lambda_{t,x} dt
\end{displaymath}
Now we shall prove that the family of measure $\{\mu_T\}_{T\geq 0}$ is tight.
So, denote with $B_R$ the set $ B_R = \{ x \in H^{1/4}: |x|_{2,1/4}^2 \leq R \} $.
Notice that since $H^{1/4} \subset H$ with compact embedding, the set $B_R$ is compact in $H$.
Moreover we have
\begin{eqnarray*}
 |X(t,x)|_{1/4}^2 - 2|W_A(t)|_{1/4}^2 &\leq& 2|Y(t,x)|_{1/4}^2 \leq \\
 &\leq& 2|2(-A+I)^{1/4}Y(t,x)|_2^2 \leq 8 |D_\xi Y(t,x)|_2^2,
\end{eqnarray*}
where $Y(t,x)$ is a strong solution of \eqref{E.00}.
Setting $\varepsilon=1$ in \eqref{E1.1}, \eqref{E1.2}, \eqref{E1.3}, it is clear that there exists a constant $C(x)$, depending by $x$, such that
\begin{displaymath}
  \int_0^T |D_\xi Y(t,x)|_2^2 dt \leq C(x)(1+T).
\end{displaymath}
Then 
\begin{eqnarray*}
  \mu(B_R^c) &=& \frac{1}{T}\int_0^T \lambda_{t,x}(B_R^c)dt = \frac{1}{T}\int_0^T \PP(|X(t,x)|_{1/4}^2 > R)dt \leq \\
  &\leq& \frac{1}{TR} \int_0^T \EE(|X(t,x)|_2) dt \leq \frac{8}{TR} \int_0^T \EE |D_\xi Y(t,x)|_2^2 dt + \\
   & & +  \frac{2}{TR} \int_0^T \EE |W_A(t)|_{1/4}^2 dt 
 \leq \frac{8}{R} C(x)\frac{(1+T)}{T}  + \frac{2}{R} \sup_{t>0} \EE|W_A(t)|_{1/4}^2.
\end{eqnarray*}
So,  it follows that $\{\mu_T \}_{t\geq 0}$ is tight.
Now, by the Krylov-Bogoliubov theorem it follows that there exists an invariant measure for the transition semigroup $P_t$. The theorem is proved.
\qed \end{proof}
%
%
\begin{lem}
  Assume that hypotesis 5.1 holds. Then for any $n\in \Nset$ we have   
\begin{equation} \label{moments}
 \int_H |x|_2^{2n} \nu(dx) < + \infty
\end{equation}
\end{lem}
\begin{proof}
 Fix $n \in \Nset$. For any $t >0$ we have
\begin{displaymath}
  |X(t,x)|_2^{2n} \leq 2n |Y(t,x)|_2^{2n} + 2n |W_A(t)|_2^{2n},
\end{displaymath}
where $Y(t,x)$ is the solution of problem \eqref{E.00}. 
Setting $\varepsilon = (1+\| DF \|_0^2/4)/2$ in \eqref{E1.1}, \eqref{E1.2}, \eqref{E1.3} it is clear that there exist
$\gamma_n, c_n > 0$ such that
\begin{displaymath}
  \EE|Y(t,x)|_2^{2n} \leq c_n (1 + \e^{-\gamma_n t }|x|_2^{2n} ).
\end{displaymath} 
Then, since $W_A(t)$ is a gaussian random variable, it follows than for some $c_n' > 0$ it holds
\begin{displaymath}
   \EE|X(t,x)|_2^{2n} \leq c_n'( 1 + \e^{-\gamma_n t }|x|_2^{2n} ).
\end{displaymath}
Now denote with $\lambda_{t,x}$ the law of $X(t,x)$. For any $\alpha > 0$ it holds
\begin{eqnarray*}
  \int_H \frac{ |y|_2^{2n} }{ 1 + \alpha |y|_2^{2n} } \lambda_{t,x}(dy) &\leq&
   \int_H  |y|_2^{2n} \lambda_{t,x}(dy) = \\
  &=& \EE|X(t,x)|_2^{2n} \leq c_n'( 1 + \e^{-\gamma_n t }|x|_2^{2n} ).
\end{eqnarray*}
Since $\lambda_{t,x} $ converges weakly to $\nu$, it follows that
\begin{displaymath}
  \int_H \frac{ |y|_2^{2n} }{ 1 + \alpha |y|_2^{2n} } \lambda_{t,x}(dy) \leq c_n'.
\end{displaymath}
Letting $\alpha \to 0 $ yields \eqref{moments}.
\qed \end{proof}
%
%
%
\section{Kolmogorov equations}
We are concerned with the semigroup $P_t$ in $L^2(H,\nu)$, where $\nu$ is the unique invariant measure for $P_t$.
In the following we only assume that $\int_H |x|_2^2 d\nu ,\infty$.
We denote by $K_2$ the infinitesimal generator of $P_t$ in  $L^2(H,\nu)$ 
and by $\mathcal{E}_A(H)$ is linear span of the set of the functions
\begin{displaymath}
  x\mapsto \cos (\langle x, h \rangle), x\mapsto \sin(\langle x, h
  \rangle),\, x\in H, h\in D(A^*).
\end{displaymath}
Let us consider the Kolmogorov operator
\begin{displaymath}
 K_0\varphi = L\varphi -  \langle F(x), D_\xi D\varphi \rangle,
 \quad \varphi \in \mathcal{E}_A(H)
\end{displaymath}
 where
\begin{displaymath}
   L\varphi(\cdot) =  \frac{1}{2}Tr[D^2\varphi(\cdot)] + \langle \cdot, AD\varphi(\cdot) \rangle
\end{displaymath}
is the Ornstein-Uhlenbek generator (see \cite{DPZ}).
Notice that $\mathcal{E}_A(H) \subset L^2(H,\nu)$ and $\mathcal{E}_A(H)$ is dense in $L^2(H,\nu)$, since $D(A)$ is dense in $H$.
Our aim is to prove that $K_2=\overline{K_0}$ in $L^2(H,\nu)$. 
First we have
\begin{lem} \label{lemma9}
For any $\varphi \in \mathcal{E}_A(H)$ we have $\varphi \in D(K_2)$ and $K_2\varphi= K_0 \varphi$.
\end{lem}
 \begin{proof}
  By It\^o's formula it follows that for all  $\varphi \in \mathcal{E}_A(H)$
\begin{displaymath}
  \lim_{t \searrow 0} \frac{1}{t} \big( P_t\varphi(x) - \varphi(x) \big) =
  K_0\varphi(x),\quad x\in H
\end{displaymath}
pointwise. Now it is sufficient to show that $\frac{1}{t} \big( P_t\varphi(x) - \varphi(x)\big)$, $ t\in (0,1]$
   is equibounded in  $L^2(H,\nu)$.
 For all $\varphi \in \mathcal{E}_A(H)$ and $x\in H$  we have
\begin{multline*}
    |P_t\varphi(x) - \varphi(x)|_2^2\leq  \\
   \leq  2t\int_0^t \EE |L\varphi(X(s,x))|_2^2 ds
      + 2t\int_0^t\EE|\langle D_\xi D\varphi(X(s,x)),F(X(s,x))\rangle |_2^2 ds.
 \end{multline*}
It is clear that there exist two positive constants $a,b$
(depending on $\varphi$) such that for all $x\in H$ it holds
\begin{displaymath}
 | L\varphi(x) |_2 \leq a+b|x|_2, \quad 
  |\langle x, D_\xi D\varphi(x)\rangle| \leq b|x|_2.
\end{displaymath}
Then we have
\begin{displaymath}
 |P_t\varphi_h(x) - \varphi_h(x)|_2^2 \leq 2 t \int_0^t \EE(a+b|X(s,x)|_2)^2 ds + 2bt \int_0^t|F(X(s,x))|_2^2ds.
\end{displaymath}
Integrating with respect to $\nu$ and taking into account the invariance of $P_t$ with respect to $\nu$ yields
\begin{displaymath}
\|P_t\varphi_h - \varphi_h\|_{L^2(H,\nu)}^2 \leq  2 t^2 \big( \int_H
(a+b|x|_2)^2 \nu (dx) + b\int_H|F(x)|_2^2 d\nu(dx) \big).
\end{displaymath}
Since $\int_H |x|_2^2 d\nu < \infty$ by assumption, the equiboundness of $t^{-1} (P_t\varphi_h -
\varphi_h)$ follows easily.
\qed \end{proof}
Before concluding that $K_2 = \overline{K_0}$ we need  two lemmas
%
%
%
%
%
\begin{lem}
  There exists a constant $c_1>0$ such that for all $h \in H$ we have 
    \begin{equation} \label{lemmaA}
      |\eta^{D_\xi h} (t,x)|_2 \leq c_1 \big(t^{-1/2} + \e^{\theta t}\big)|h|_2,
    \end{equation}
where $\theta = \big( \kappa \|DF\|_0 \Gamma(1/2) \big)^2$ and $\eta^z(t,x)$ is defined as in \eqref{ED.1}.
\end{lem}
\begin{proof}
 Notice that by the density of $H^1$ in $H$ it is sufficient to prove \eqref{lemmaA} for $h\in H^1$.
 So, if $h\in H^1$,
 $\eta^{D_\xi h}(t,x)$ is the solution of
\begin{equation*}
 \eta^{D_\xi h}(t,x) = K(t) h +
     \int_0^t K(t-s)\langle DF(X(s,x)),\eta^{D_\xi h}(s,x)\rangle
     ds.
\end{equation*}
By \eqref{2.3} it follows that
\begin{displaymath}
 |\eta^{D_\xi h}(t,x)|_2 \leq \kappa t^{-1/2}|h|_2 +
 \kappa \|DF\|_0 \int_0^t  (t-s)^{-1/2}  |\eta^{D_\xi h}(s,x)|_2 ds
\end{displaymath}
and by Lemma \ref{Henry} that 
\begin{equation} \label{lemmaA1}
 |\eta^{D_\xi h}(t,x)|_2 \leq  \kappa t^{-1/2}|h|_2 +
 \kappa \theta \int_0^t E'_{1/2}(\theta(t-s))s^{-1/2}ds |h|_2.
\end{equation}
Since $E_{1/2}(\cdot):\Rset^+ \to \Rset$ is a continuous function with
\begin{displaymath}
 E'_{1/2}(z) \sim  2z^{-1/2} \quad \text{as $z \to 0^+$ },\quad E'_{1/2}(z) \sim 2 \e^z \quad \text{as $z\to +\infty$},
\end{displaymath}
it is clear that there exists a constant $C>0$ such that 
\begin{displaymath}
  E'_{1/2}(z) \leq C(z^{-1/2} + \e^z) \quad z>0.
\end{displaymath}
Taking into account that 
\begin{displaymath}
  \int_0^t (t-s)^{-1/2}s^{-1/2}ds =
    \int_0^1 (1-s)^{-1/2}s^{-1/2}ds =
    \beta(1/2,1/2),
\end{displaymath}
where $\beta(\cdot,\cdot)$ is the Euler beta funtion, by an easy computation we find \eqref{lemmaA}
\qed \end{proof}
%
%
%
%
\begin{lem}
 Assume that $f \in C^1_b(H)$. Then there exists $c_2 > 0$ such
that for all $\lambda> \theta $, $h \in H$ we have 
  \begin{equation} \label{lemmaB}
   |\langle D_\xi D \varphi(x), h \rangle|  \leq
     c_2 \big(\lambda^{-1/2} + \frac{1}{\lambda - \theta} \big) \|f\|_1 |h|_2
  \end{equation}
 where
  \begin{displaymath}
      \theta = \big(\kappa \|DF\|_0 \Gamma(1/2)\big)^2
  \end{displaymath}
  and
  \begin{displaymath}
    \varphi(x) = \int_0^\infty \e^{-\lambda t} P_t f(x) dt.
  \end{displaymath}
\end{lem}
\begin{proof}
 Notice that since $H^1$ is dense in $H$ it is sufficient to prove \eqref{lemmaB} for $h\in H^1$. We have
 \begin{eqnarray*}
    \langle D_\xi D \varphi(x), h \rangle &=& - \langle  D \varphi(x),
    D_\xi h \rangle = \\
  &=& - \int_0^\infty \e^{-\lambda t} \EE \big[ \langle Df(X(t,x)), \langle
   X_x(t,x),D_\xi h \rangle \rangle\big] dt = \\
  &=& -\int_0^\infty \e^{-\lambda t} \EE \big[ \langle Df(X(t,x)),
   \eta^{D_\xi h} (t,x)  \rangle\big]  dt.
  \end{eqnarray*}
 So, taking into account \eqref{lemmaA},
\begin{eqnarray*}
  | \langle D_\xi D \varphi(x), h \rangle|_2  &\leq& \|f\|_1  \int_0^\infty \e^{-\lambda t}
   |\eta^{D_\xi h} (t,x)|_2   dt \leq \\
  &\leq&  \|f\|_1  c_1 \int_0^\infty \e^{-\lambda t} (t^{-1/2} + \e^{\theta t } )dt |h|_2
 \end{eqnarray*}
 that implies \eqref{lemmaB}.
\qed \end{proof}
%
%
%
%
Now we are able to prove the main result of this section.
\begin{thm}
 $K_2$ is the closure of $K_0$ in $L^2(H,\nu)$.
\end{thm}
\begin{proof}
 By Lemma \ref{lemma9} we know that $K_2$ extends $K_0$. Since $K_2$ is dissipative, so is $K_0$. 
 Consequently, $K_0$ is closable. Let us denote by $\overline{K_0}$ its closure. 
 We have to show that $K_2 = \overline{K_0}$. 
 For this pourpose, we will show that the
 range of $\lambda - \overline{K_0}$ is dense in $L^2(H,\nu)$ for some $\lambda > 0$. 
 In fact by the Lumer-Philipps theorem this implies that $\overline{K_0}$ is $m$-dissipative and it is the generator of a semigroup of contraction. 
 Therefore, since $K_2$ extends $\overline{K_0}$, it must coincide with $\overline{K_0}$. 
 So, if $f\in \mathcal{E}_A(H)$, we are interested in solving the problem
\begin{displaymath}
 \lambda \varphi - \overline{K_0} \varphi = f.
\end{displaymath}
Setting
\begin{displaymath}
   \varphi (x) = R(\lambda,K_2)f(x) = \int_0^\infty \e^{-\lambda t} P_t f(x) dt,
\end{displaymath}
 we will show in the next two steps that $\varphi \in D(\overline{K_0})$, that implies
 $\lambda \varphi - \overline{K_0} \varphi = f$, since $K_2$ extends $\overline{K_0}$.
 Let us denote with $C_{b,1}(H)$ the Banach space of all continuous function $\psi:H\to \Rset$ such that $\psi(x)/(1+|x|_2)$ $\in$ $C_b(H)$.\\
\textbf{Step 1.} $\varphi \in D(L,C_{b,1}(H))$ \\
 Notice that by Theorem \ref{SF} it follows that $\varphi \in C_b^1(H)$. We have to compute the derivative
\begin{displaymath}
  \frac{d}{dt} R_t\varphi|_{t=0},
\end{displaymath}
 where
\begin{displaymath}
   R_t\varphi (x) = \EE\big[\varphi(Z(t,x))\big], \quad
   Z(t,x) = \e^{tA}x + W_A(t).
\end{displaymath}
We have  
\begin{multline} \label{step1}
  R_t\varphi (x) = \EE\big[\varphi(Z(t,x))\big] =
  \EE\big[\varphi\big(X(t,x) - \int_0^t K(t-s)F(X(s,x))ds\big)  \big] =   \\
  = P_t\varphi(x) - \EE \big[ \langle D\varphi(X(t,x)), \int_0^t
  K(t-s)F(X(s,x))ds \rangle \big] + o(t),
\end{multline}
where $\lim_{t\to 0} \frac{o(t)}{t} = 0 $.
Now, since $\forall x \in H$ we have that $\PP-$a.s.
\begin{equation} \label{step1.0}
   \lim_{t\to 0^+} \frac{1}{t}  \langle D\varphi(X(t,x)), \int_0^t
  K(t-s)F(X(s,x))ds \rangle = - \langle D_\xi D\varphi(x), F(x) \rangle 
\end{equation}
Taking into account \eqref{lemmaB} with $\lambda > \theta$ we find
\begin{eqnarray}\label{step1.1}
  \frac{1}{t}  |\langle D\varphi(X(t,x)) &,& \int_0^t   K(t-s)F(X(s,x))ds \rangle | = \notag \\
 &=& \e^{\lambda t} \frac{1}{t}  |\langle D\varphi(x), \int_0^t   K(t-s)F(X(s,x))ds \rangle | =\notag\\
 &=& \e^{\lambda t} \frac{1}{t}  |\langle D_\xi  D\varphi(x), \int_0^t  \e^{(t-s)A} F(X(s,x))ds \rangle | \leq\notag\\
 &\leq& \e^{\lambda t} c_2(\lambda^{ 1/2} + \frac{1}{\lambda- \theta} ) \| DF \|_0 \frac{1 - \e^{-t}}{t} \sup_{0\leq s\leq T} |X(s,x)|_2
\end{eqnarray}
So, since \eqref{step1.0}, \eqref{step1.1} hold, we can apply the dominated convergence theorem in \eqref{step1} and obtain
\begin{displaymath}
  \lim_{t\to 0} \frac{R_t \varphi(x)- \varphi(x)}{t} = 
   \lambda \varphi(x) - f(x) +  \langle D_\xi D\varphi(x), F(x) \rangle
\end{displaymath}
 for all $x \in H$.
 Now we have to prove that $ t^{-1}\big(R_t \varphi(x) - \varphi(x)\big) $
 is equibounded in $C_{b,1}(H)$ for $t\in (0,1]$. 
 To see this we need to observe that by \eqref{2.1} it follows easily that 
\begin{displaymath} 
  \sup_{t\in [0,1]} \EE|X(t,x)|_2 \leq |x|_2 C, 
\end{displaymath}  
where $C = \sup_{t \geq 0 } \EE|W_A(t)|_2 E_{1/2}(\theta)$. Now set 
\begin{displaymath}
C' = c_2 (\lambda^{-1/2} + \frac{1}{\lambda - \theta})\|DF\|_0 C.
\end{displaymath}
By \eqref{step1} and \eqref{step1.1} we have
\begin{multline*}
  \frac{|R_t\varphi(x) - \varphi(x)|}{t(1+|x|_2)} \leq  \frac{|P_t\varphi(x) - \varphi(x)|}{t(1+|x|_2)} + \\
  - \frac{1}{t(1+|x|_2)}\EE\big[\langle D\varphi(X(t,x)), \int_0^t
  K(t-s)F(X(s,x))ds \rangle] + \frac{o(t)}{t} \leq  \\
  \leq \frac{\e^{\lambda t} - 1}{t} \|\varphi\|_{b,1} + C'\frac{\e^{\lambda t}(1-\e^{-t})|x|_2}{t(1+|x|_2)} +\frac{o(t)}{t}
\end{multline*}
which is equibounded in $(0,1]$. \\
\textbf{Step 2.} $D(L,C_{b,1}(H)) \subset D(\overline{K_0})$.\\
 By Proposition 2.6 of \cite{DPT} there exists a 4-index sequence
  $(\varphi_{n_1,n_2,n_3,n_4}) \subset \mathcal{E}_A(H)$ such that for all $x \in H$
\begin{align*}
 \lim_{n_1 \to \infty} \lim_{n_2 \to \infty}\lim_{n_3 \to \infty}\lim_{n_4 \to \infty}\varphi_{n_1,n_2,n_3,n_4} (x) &=
 \varphi(x) \\
 \lim_{n_1 \to \infty} \lim_{n_2 \to \infty}\lim_{n_3 \to \infty}\lim_{n_4 \to \infty} L\varphi_{n_1,n_2,n_3,n_4} (x) &=
    L\varphi(x) \\
 \lim_{n_1 \to \infty} \lim_{n_2 \to \infty}\lim_{n_3 \to \infty}\lim_{n_4 \to \infty} D\varphi_{n_1,n_2,n_3,n_4} (x) &=
    D\varphi(x)
\end{align*}
and
\begin{displaymath}
  \sup_{n_1,n_2,n_3,n_4} \{ \|\varphi_{n_1,\,n_2,\,n_3,n_4}\|_{b,2} +
                         \|D\varphi_{n_1,\,n_2,\,n_3,n_4}\|_{b,2} +
                       \|L\varphi_{n_1,\,n_2,\,n_3,n_4}\|_{b,2} \} < \infty.
\end{displaymath}
So, by the dominated convergence theorem it follows that
\begin{displaymath}
 \lim_{n\to\infty} K_0 \varphi_n =\lim_{n \to\infty}( L \varphi_n -
 \langle D_\xi D \varphi_n , F \rangle)
 = L\varphi - \langle D_\xi D\varphi, F \rangle
\end{displaymath}
in $L^2(H,\nu)$, that implies $\varphi \in
D(\overline{K_0})$.
The theorem is proved.
\qed \end{proof}
%
%
%
%
\section{Integration by parts formula}         
Let us denote by $\nu$ the invariant measure for the transition
semigroup $P_t$ and by $K_2$ its infinitesimal generator in
$L^2(H,\nu)$.
\begin{prop}
 The operator $D:\mathcal{E}_A(H) \to C_b(H,H)$, $\varphi \mapsto
 D\varphi$, is uniquely extendible to a linear bounded operator
 $D:D(K_2) \to L^2(H,\nu;H)$.
 Moreover, \eqref{IBPF} holds.
\end{prop}
\begin{proof}
 Let $\varphi \in \mathcal{E}_A(H)$.
 Taking into account Lemma \ref{lemma9} and that $\varphi^2 \in \mathcal{E}_A(H)$ by a simple computation we see that
$ K_2(\varphi^2) = 2 \varphi K_2 \varphi + |D\varphi|_2^2$.
Integrating both sides over $H$ with respect to $\nu$ and taking into account that
\begin{displaymath}
  \int_H K_2(\varphi^2) d\nu = 0
\end{displaymath}
by the invariance of $P_t$ with respect to $\nu$, it follows \eqref{IBPF}. Now we will prove that \eqref{IBPF} holds for all $\varphi \in D(K_2)$. Let us fix $\varphi \in D(K_2)$.
 Since $\mathcal{E}_A(H)$ is a core for $P_t$, there exists a sequence $\{\varphi_n\} \subset \mathcal{E}_A(H)$
such that
\begin{displaymath}
  \varphi_n \to \varphi,\; K_2\varphi_n \to K_2 \varphi \quad \text{in $L^2(H,\nu)$},
\end{displaymath}
and consequently
\begin{eqnarray*}
 \int_H |D(\varphi_n -\varphi_m)|_2^2 d\nu &\leq& 
   2 \int_H |\varphi_n -\varphi_m|_2|K_2(\varphi_n -\varphi_m)|_2d\nu \leq\\
     &\leq& \| \varphi_n -\varphi_m \|_{L^2(H;\nu)}^2 + \|K_2(\varphi_n -\varphi_m)\|_{L^2(H;\nu)}^2.
\end{eqnarray*}
Therefore the sequence $\{ D\varphi_n \}$ is Cauchy in $L^2(H,\nu;H)$, and the conclusion follows.
\qed \end{proof}
\section{The Sobolev space $W^{1,2}(H,\nu)$}

 We want to show that the mapping
\begin{displaymath}
D:\mathcal{E}_A(H)\subset L^2(H,\nu)\to L^2(H,\nu;H),\; \varphi\mapsto D\varphi
\end{displaymath}
is closable.  
\begin{thm} 
 \label{t3.24}
  $D$ is closable. Moreover, if $\varphi$ belongs to the domain $\overline{D}$ of the closure of $D$ and 
 $\overline{D}\varphi=0$ we have that $\overline{D}P_t\varphi=0$ for any $t>0.$
\end{thm}
\begin{proof}
Since $\|R(\lambda,L)\|_{\mathcal{L}(L^2(H;\nu))} \leq \sqrt{\lambda/t},\, \lambda > 0$ and \eqref{ED.1} holds, we can apply Proposition 3.5 of \cite{DPDG}. The theorem is proved.
\qed \end{proof}
We define by $W^{1,2}(H,\nu)$ the domain of $D$ in $L^2(H,\nu)$.
   By \eqref{IBPF} it follows that $D(K_2) \subset  W^{1,2}(H,\nu)$.
   We have the next
\begin{prop}
Let $\varphi \in L^2(H,\nu)$ and $t\geq 0$. Set $u(t,x)=P_t\varphi(x)$. Then, for any $T >0$,
we have $u \in L^2(0,T;W^{1,2}(H,\nu))$ and the following indentity holds 
\begin{equation} \label{PF2}
  \int_H (P_t \varphi)^2 d\nu + \int_0^t ds \int_H |DP_s \varphi|^2 d\nu = \int_H \varphi^2 d\nu.
\end{equation}
\end{prop}
\begin{proof}
 Let first $\varphi \in D(K_2)$. We have that $P_t \varphi \in D(K_2)$ and for all $t\geq 0$
 \begin{displaymath}
    \frac{d}{dt} P_t\varphi(x) = K_2 P_t\varphi(x).
 \end{displaymath}
Multiplying both sides of this identity by $P_t\varphi(x)$ and integrating with respect to $x$ over $H$,
   by \eqref{IBPF}  yields
\begin{displaymath}
 \frac{1}{2} \frac{d}{dt} \int_H (P_t \varphi)^2 d\nu = \int_H P_t \varphi K_2 P_t\varphi d\nu =
 -\frac{1}{2} \int_H |D P_t \varphi |^2 d\nu.
\end{displaymath}
Integrating with respect to $t$ it yields \eqref{PF2}.
  Now the conclusion follows by the density of $D(K_2)$ in $L^2(H,\nu)$.
\qed \end{proof}
Letting $t\to \infty$ in \eqref{PF2} we have
\begin{prop}
 For any $\varphi \in L^2(H,\nu)$, we have 
\begin{equation} \label{PF3}
 \int_H |\varphi - \overline{\varphi}|^2 d\nu = \int_0^\infty dt \int_H |DP_t\varphi|^2d\nu,
\end{equation}
where $\overline{\varphi} = \int_H \varphi d\nu$.
\end{prop}
%
%
\section{Poincar\'e and log-Sobolev inequality, spectral gap}
As in  \eqref{PF3} we will use the notation $   \overline{\varphi} = \int_H \varphi d\nu.$
Let us prove the Poincar\'e inequality.
\begin{thm}[Poincar\'e inequality]
 Let us assume  $\|DF\|_0 < 2$.
 Then for any $\varphi \in W^{1,2}(H,\nu)$ inequality \eqref{PI} holds. 
\end{thm}
\begin{proof}
 Let first $\varphi \in \mathcal{E}_A(H)$.
 Then for any $h\in H$ and $t\geq 0$ we have
\begin{displaymath}
  \langle DP_t\varphi (x), h\rangle = \EE [\langle D\varphi(X(t,x)),\eta^h(t,x)\rangle ]
\end{displaymath}
 Taking into account \eqref{ED.3} it follows that
\begin{multline*}
   |\langle DP_t\varphi (x), h\rangle|^2 \leq
     \EE [ |D\varphi(X(t,x))|_2^2 |\eta^h(t,x)|_2^2 ] \leq \\
   \leq \e^{2\big(\frac{\|DF\|_0^2}{4} - 1 \big)t } \EE[|D\varphi (X(t,x))|_2^2] |h|_2^2 =
        \e^{2\big(\frac{\|DF\|_0^2}{4} - 1 \big)t }P_t(|D\varphi|_2^2)(x) |h|_2^2.
\end{multline*}
By the arbitrariness of $h$ it yields
\begin{displaymath}
 |DP_t\varphi (x)|_2^2 \leq \e^{2\big(\frac{\|DF\|_0^2}{4} - 1 \big)t }P_t(|D\varphi|_2^2)(x)
\end{displaymath}
 for all $x\in H$, $s\geq 0$.
Taking into account \eqref{PF3} and the invariance of $\nu$, we obtain
\begin{eqnarray*}
  \int_H |\varphi - \overline{\varphi}|^2 d\nu
   &\leq& \int_0^\infty dt \int_H \e^{2\big(\frac{\|DF\|_0^2}{4} - 1 \big)t }P_t(|D\varphi|_2^2) d\nu = \\
   &=& \frac{1}{2\big( 1- \frac{\|DF\|_0^2}{4} \big)} \int_H |D\varphi|_2^2 d\nu
\end{eqnarray*}
and the conclusion follows. If $\varphi \in W^{1,2}(H,\nu)$ we proceed by density.
\qed \end{proof}
%
%
%

Now we show that if \eqref{PI} holds then there is a gap in the spectrum of $K_2$ and that the convergence to the equilibrium point is exponential. We have
\begin{thm}[Spectral gap]
 Let us assume  $\|DF\|_0 <2$. Then we have 
\begin{equation}\label{SG1}
  \sigma(K_2) \setminus  \{  0  \} \subset
   \{ \lambda \in \Cset: \Re e\lambda \leq - \big( 1- \frac{\|DF\|_0^2}{4} \big) \}.
\end{equation}
Moreover     
\begin{equation} \label{SG}
   \int_H |P_t\varphi - \overline{\varphi}|^2 d\nu \leq
   \e^{-2\big( 1- \frac{\|DF\|_0^2}{4} \big)t}\int_H |\varphi|^2 d\nu.
\end{equation}
\end{thm}
\begin{proof}
The proof is an easy consequence of \eqref{IBPF} and \eqref{PI} ( see \cite{DPDG} and \cite[Prop. 2.3]{GR}).
\qed \end{proof}
%
%
\begin{thm}[Log-Sobolev inequality]
 Let us assume  $\|DF\|_0 <2$.
Then or any $\varphi \in W^{1,2}(H,\nu)$ we have 
\begin{equation} \label{logsob} 
  \|\varphi^2 \log(\varphi^2)\|_{L^2(H,\nu)} \leq \frac{1}{1-\frac{\| DF \|_0^2}{4}}  \|D\varphi\|_{L^2(H,\nu)}^2 + 
  \|\varphi^2\|_{L^2(H,\nu)} \log(\|\varphi^2\|_{L^2(H,\nu)}).
\end{equation}
Moreover, the transition semigroup $\{P_t\}_{t\geq 0}$ is hypercontractive.
\end{thm}
\begin{proof}
Let us take $\varphi\in \mathcal{E}_A(H)$, with $\varphi \neq 0$. 
We have
\begin{displaymath}
  \frac{d}{dt} \int_H P_t(\varphi^2) \log(P_t(\varphi^2)) d\nu = \int_H K_2 P_t(\varphi^2) \log(P_t(\varphi^2))d\nu +
     \int_H K_2 P_t(\varphi^2)d\nu,
\end{displaymath}
where the last term  vanishes  due to the invariance of $\nu$. Moreover, it holds the identity 
\begin{equation} \label{logsob1}
  \int_H g'(\varphi) K_2 \varphi  d\nu = -\frac{1}{2}  \int_H g''(\varphi) |D\varphi|^2 d\nu.
\end{equation}
Since 
\begin{displaymath}
    D P_t(\varphi^2) = 2 \EE[\varphi(X(t,x)) D\varphi(X(t,x))\cdot X_x(t,x)]
\end{displaymath}
it follows, from the H\"older inequality and \eqref{ED.4} that
\begin{eqnarray*}
 |D P_t(\varphi^2)|^2 &\leq& 4 \EE[|\varphi|^2(X(t,x))]\EE[|D\varphi|^2(X(t,x))] \e^{-2(1-\|DF\|_0^2/4)t}=\\
 &=& 4 \e^{-2(1-\|DF\|_0^2/4)t} P_t(|\varphi|^2)(x) P_t(|D\varphi|^2)(x).
\end{eqnarray*}
Therefore, by \eqref{logsob1} it yields
\begin{eqnarray*}
   \frac{d}{dt} \int_H P_t(\varphi^2) \log(P_t(\varphi^2)) d\nu &\geq& 
    -2 \e^{-2(1-\|DF\|_0^2/4)t} \int_H P_t(|D\varphi|^2)(x) d\nu = \\
   &=&  -2 \e^{-2(1-\|DF\|_0^2/4)t} \int_H |D\varphi|^2(x) d\nu 
\end{eqnarray*}
due to the invariance of $\nu$.
Integrating with respect to $t$ gives
\begin{multline*}
  \int_H P_t(\varphi^2) \log(P_t(\varphi^2)) d\nu   \geq  \\
  \geq \int_H \varphi^2 \log(\varphi^2) d\nu  - 
   \frac{1- \e^{-2(1-\|DF\|_0^2/4)t}}{ 2(1-\|DF\|_0^2/4) } \int_H |D\varphi|^2(x) d\nu. 
\end{multline*}
Then the conclusion follows letting $t \to \infty $.
Finally, as shown in \cite{GR}, a logarithmic Sobolev inequality implies that the transition semigroup $\{P_t\}_{t\geq 0}$ is hypercontractive.

\qed \end{proof}


\begin{thebibliography}{5}

\bibitem{CER} Cerrai, S.,  \emph{``Second order PDE's in finite and infinite dimensions. A probabilistic approach''}. Lecture Notes in Mathematics, 1762, Springer-Verlag, 2001.
\bibitem{DP03} Da Prato, G., \emph{``Kolmogorov Equations for Stochastic PDEs''}. Advanced Courses in Mathematics. CRM Barcelona. Birkhäuser Verlag, Basel, 2004.
\bibitem{DPD} Da Prato, G.; Debussche, A.,  Maximal dissipativity of the Dirichlet operator corresponding to the Burgers equation.  \emph{Stochastic processes, physics and geometry: new interplays, I (Leipzig, 1999)},  85--98, CMS Conf. Proc., 28, Amer. Math. Soc., Providence, RI, 2000.
\bibitem{DPEZ} Da Prato, G.; Elworthy, D.; Zabczyk, J., Strong Feller property for stochastic semilinear equations. Stoch. Anal. Appl., \textbf{1995}, 13, 35-45
\bibitem{DPZ96} Da Prato, G.; Zabczyk, J.,  \emph{``Ergodicity for Infinite Dimensional System''}. London Mathematical Society Lecture Notes Series, 229, Cambridge University Press, Cambridge, 1996
\bibitem{DPZ02} Da Prato, G.; Zabczyk, J., \emph{``Second Order Partial Differential Equations in Hilbert Spaces''}. London Mathematical Society Lecture Notes, 293, Cambridge University Press, Cambridge, 2002.
\bibitem{DPZ} Da Prato, G.; Zabczyk, J.,  \emph{``Stochastic equations in infinite dimensions''}. Cambridge University Press, Cambridge, 1992
\bibitem{DPT} Da Prato, G.; Tubaro, L., Some results about dissipativity of Kolmogorov operators. Czechoslovak Math. J., \textbf{2001}, 51(126),  685-699.
\bibitem{DPDG} Da Prato, G.; Debussche, A.; Goldys, B., Invariant measures of non symmetric dissipative stochastic systems. Probab. Theory Relat. Fields, \textbf{2002}, 123, no. 3, 355-380.
\bibitem{GR} Gross, L., Logarithmic Sobolev inequalities. Amer. J. Math., \textbf{1975}, 97, no. 4, 1061-1083
\bibitem{GY98} Gy\"ongy, I., Existence and uniqueness results for semilinear stochastic partial differential equations. Stochastic Process. Appl., , \textbf{1998}, 78, no. 2, 271-299
\bibitem{Hen} Henry, D., \emph{``Geometric theory of semilinear parabolic  equations''}. Springer-Verlag Lecture Notes in Mathematics, 840, Berlin-New York, 1981
\end{thebibliography}
\end{document}